\documentclass[12pt]{article}
\usepackage{amssymb}
\usepackage{amsfonts}
\usepackage{latexsym}
\usepackage{amsthm}

\newcommand{\er}[1]{{\rm(\ref{#1})}}
\def\lb{\label}

\theoremstyle{plain}
\newtheorem{theorem}{\bf Theorem}[section]
\newtheorem{lemma}[theorem]{\bf Lemma}

\newtheorem{proposition}[theorem]{\bf Proposition}
\theoremstyle{remark}

\begin{document}

\def\a{\alpha}  \def\cA{{\cal A}}     \def\bA{{\bf A}}  \def\mA{{\mathscr A}}
\def\b{\beta}   \def\cB{{\cal B}}     \def\bB{{\bf B}}  \def\mB{{\mathscr B}}
\def\g{\gamma}  \def\cC{{\cal C}}     \def\bC{{\bf C}}  \def\mC{{\mathscr C}}
\def\G{\Gamma}  \def\cD{{\cal D}}     \def\bD{{\bf D}}  \def\mD{{\mathscr D}}
\def\d{\delta}  \def\cE{{\cal E}}     \def\bE{{\bf E}}  \def\mE{{\mathscr E}}
\def\D{\Delta}  \def\cF{{\cal F}}     \def\bF{{\bf F}}  \def\mF{{\mathscr F}}
\def\c{\chi}    \def\cG{{\cal G}}     \def\bG{{\bf G}}  \def\mG{{\mathscr G}}
\def\z{\zeta}   \def\cH{{\cal H}}     \def\bH{{\bf H}}  \def\mH{{\mathscr H}}
\def\e{\eta}    \def\cI{{\cal I}}     \def\bI{{\bf I}}  \def\mI{{\mathscr I}}
\def\p{\psi}    \def\cJ{{\cal J}}     \def\bJ{{\bf J}}  \def\mJ{{\mathscr J}}
\def\vT{\Theta} \def\cK{{\cal K}}     \def\bK{{\bf K}}  \def\mK{{\mathscr K}}
\def\k{\kappa}  \def\cL{{\cal L}}     \def\bL{{\bf L}}  \def\mL{{\mathscr L}}
\def\l{\lambda} \def\cM{{\cal M}}     \def\bM{{\bf M}}  \def\mM{{\mathscr M}}
\def\L{\Lambda} \def\cN{{\cal N}}     \def\bN{{\bf N}}  \def\mN{{\mathscr N}}
\def\m{\mu}     \def\cO{{\cal O}}     \def\bO{{\bf O}}  \def\mO{{\mathscr O}}
\def\n{\nu}     \def\cP{{\cal P}}     \def\bP{{\bf P}}  \def\mP{{\mathscr P}}
\def\r{\rho}    \def\cQ{{\cal Q}}     \def\bQ{{\bf Q}}  \def\mQ{{\mathscr Q}}
\def\s{\sigma}  \def\cR{{\cal R}}     \def\bR{{\bf R}}  \def\mR{{\mathscr R}}
\def\S{\Sigma}  \def\cS{{\cal S}}     \def\bS{{\bf S}}  \def\mS{{\mathscr S}}
\def\t{\tau}    \def\cT{{\cal T}}     \def\bT{{\bf T}}  \def\mT{{\mathscr T}}
\def\f{\phi}    \def\cU{{\cal U}}     \def\bU{{\bf U}}  \def\mU{{\mathscr U}}
\def\F{\Phi}    \def\cV{{\cal V}}     \def\bV{{\bf V}}  \def\mV{{\mathscr V}}
\def\P{\Psi}    \def\cW{{\cal W}}     \def\bW{{\bf W}}  \def\mW{{\mathscr W}}
\def\o{\omega}  \def\cX{{\cal X}}     \def\bX{{\bf X}}  \def\mX{{\mathscr X}}
\def\x{\xi}     \def\cY{{\cal Y}}     \def\bY{{\bf Y}}  \def\mY{{\mathscr Y}}
\def\X{\Xi}     \def\cZ{{\cal Z}}     \def\bZ{{\bf Z}}  \def\mZ{{\mathscr Z}}
\def\O{\Omega}
\def\ve{\varepsilon}
\def\vt{\vartheta}
\def\vp{\varphi}

\def\Z{{\Bbb Z}}
\def\R{{\Bbb R}}
\def\C{{\Bbb C}}
\def\T{{\Bbb T}}
\def\N{{\Bbb N}}
\def\dD{{\Bbb D}}

\def\sign{\mathop{\rm sign}\nolimits}
\def\sign{\mathop{\rm sign}\nolimits}
\def\dist{\mathop{\rm dist}\nolimits}
\def\Ker{\mathop{\rm Ker}\nolimits}
\def\bs{\backslash}
\def\iint{\int\!\!\!\int}
\def\bl{\biggl}   \def\br{\biggr}

\def\ma{\left(\begin{array}{cc}}
\def\am{\end{array}\right)}
\def\intl{\int\limits}
\def\iintl{\iint\limits}
\def\tminf{{\underset{m\to\infty}\to \longrightarrow}}

\def\qqq{\qquad}
\def\qq{\quad}
\let\ge\geqslant
\let\le\leqslant
\let\geq\geqslant
\let\leq\leqslant

\def\ma{\left(\begin{array}{cc}}    \def\am{\end{array}\right)}
\def\iint{\int\!\!\!\int}
\def\lt{\biggl}                     \def\rt{\biggr}
\let\ge\geqslant                   \let\le\leqslant
\def\[{\begin{equation}}            \def\]{\end{equation}}
\def\wt{\widetilde}                 \def\pa{\partial}
\def\sm{\setminus}                  \def\es{\emptyset}
\def\no{\noindent}                  \def\ol{\overline}
\def\iy{\infty}                     \def\ev{\equiv}
\def\/{\over}
\def\ts{\times}
\def\os{\oplus}
\def\ss{\subset}
\def\h{\hat}
\def\wh{\widehat}
\def\Ra{\Rightarrow}
\def\ra{\rightarrow}
\def\la{\leftarrow}
\def\da{\downarrow}
\def\ua{\uparrow}
\def\lra{\leftrightarrow}
\def\Lra{\Leftrightarrow}
\def\Re{\mathop{\rm Re}\nolimits}
\def\Im{\mathop{\rm Im}\nolimits}
\def\supp{\mathop{\rm supp}\nolimits}
\def\sign{\mathop{\rm sign}\nolimits}
\def\Ran{\mathop{\rm Ran}\nolimits}
\def\Ker{\mathop{\rm Ker}\nolimits}
\def\Tr{\mathop{\rm Tr}\nolimits}
\def\const{\mathop{\rm const}\nolimits}
\def\Wr{\mathop{\rm Wr}\nolimits}

\def\th{\theta}
\def\dlint{\displaystyle\int\limits}
\def\iintt{\mathop{\int\!\!\int\!\!\dots\!\!\int}\limits}
\def\intt{\mathop{\int\int}\limits}
\def\lim{\mathop{\rm lim}\limits}
\def\mult{\!\cdot\!}
\def\BBox{\hspace{1mm}\vrule height6pt width5.5pt depth0pt \hspace{6pt}}
\def\1{1\!\!1}
\newcommand{\bwt}[1]{{\mathop{#1}\limits^{{}_{\,\bf{\sim}}}}\vphantom{#1}}
\newcommand{\bhat}[1]{{\mathop{#1}\limits^{{}_{\,\bf{\wedge}}}}\vphantom{#1}}
\newcommand{\bcheck}[1]{{\mathop{#1}\limits^{{}_{\,\bf{\vee}}}}\vphantom{#1}}
\def\nh{\bhat}
\def\nc{\bcheck}
\newcommand{\oo}[1]{{\mathop{#1}\limits^{\,\circ}}\vphantom{#1}}
\newcommand{\po}[1]{{\mathop{#1}\limits^{\phantom{\circ}}}\vphantom{#1}}
\def\ctg{\mathop{\rm ctg}\nolimits}
\def\notto{\to\!\!\!\!\!\!\!/\,\,\,}

\def\pgbrk{\pagebreak}

\title {Estimates for periodic Zakharov-Shabat operators}

\author{
Evgeny Korotyaev
\begin{footnote}
{Cardiff School of Mathematics, Cardiff University, Senghennydd
Road, Cardiff CF24 4AG, UK, e-mail: korotyaeve@cardiff.ac.uk}
\end{footnote}
and  Pavel Kargaev
\begin{footnote}
{Faculty of Math. and Mech. St-Petersburg State University, e-mail:
kargaev@PK2673.spb.edu}
\end{footnote}
}
\maketitle

\begin{abstract}
\no We consider the periodic Zakharov-Shabat operators on the real
line. The spectrum of this operator consists of intervals separated
by gaps with the lengths $|g_n|\ge 0, n\in \Z$. Let $\m_n^\pm$ be
the corresponding effective masses and let $h_n$ be heights of the
corresponding slits in the quasimomentum domain. We obtain a priori
estimates of sequences $g=(|g_n|)_{n\in \Z},\m^\pm=(\m_n^\pm)_{n\in
\Z},  h=(h_n)_{n\in \Z}$  in terms of weighted $\ell^p-$norms at
$p\ge 1$.  The proof is based on  the analysis of the quasimomentum
as the conformal mapping.
\end{abstract}


\section{Introduction and  main results}
\setcounter{equation}{0}

Consider the Zakharov-Shabat operator $\cT$ acting on $L^2(\R )\os
L^2(\R )$ and given by
$$
\cT=\cJ {d\/ dt}+V(t),\ \ \ \ V=\ma V_1&V_2\\V_2&-V_1\am \ , \ \
\cJ= \ma 0&1\\-1&0\am ,
 $$
where $V$ is a real 1-periodic $2\ts 2$ matrix-valued function of
$t\in \R$ and  $V\in L^1(0,1)$. In order to describe our main result
we shall introduce some notations and recall some well known facts
about the Zakharov-Shabat operator (see [Kr], [LS] for details). The
spectrum of $\cT$ is purely absolutely continuous and is given by
the set $\cup \s_n$, with spectral bands $\s_n=[z_{n-1}^+,z_n^-]$,
where $\dots <z_{2n-1}^-\le z_{2n-1}^+< z_{2n}^- \le z_{2n}^+< \dots
$ and $\  z_n^{\pm}=n(\pi +o(1))$   as $|n|\to\iy$. These intervals
$\s_n, \s_{n+1}$ are separated by a gap $g_n=(z^-_n,z^+_n)$ with
length $|g_n|\geq 0$.  If a gap $g_n$ is degenerate, i.e. $g_n=\es$,
then the corresponding spectral bands $\s_n, \s_{n+1}$ merge.  The
sequence $\dots <z_{2n-1}^-\le z_{2n-1}^+< z_{2n}^- \le z_{2n}^+<
\dots$  is the spectrum of equation $\cJ f'+Vf=zf$ with the
2-periodic boundary conditions, i.e., $f(t+2)=f(t), t\in \R$. Here
the equality $z_n^-=z_n^+$ means that $z_n^-$ is the double
eigenvalue.  The eigenfunctions, corresponding to the eigenvalue
$z_n^{\pm}$, are 1-periodic,  when $n$ is even and they are
antiperiodic,  i.e., $f(t+1)=-f(t),\ t\in \R$, when $n$ is odd.
Introduce the  $2\ts 2$-matrix valued fundamental solution
 $\p =\p (t,z)$ of the problem:
\[
\cJ\p '+V\p=z \p, \ \  \p (0,z)=\cI_2, \ \ z\in \C,
\]
where $\cI_2$ is the identity $2\times 2$-matrix. Introduce the
Lyapunov function $\D(z)={1 \/ 2}{\rm Tr} \p (1,z)$. Note that
$\D(z_n^{\pm})=(-1)^n, n\in\Z$, and the function $\D'(z)$ has
exactly one zero $z_n\in [z^-_n,z^+_n]$  for each $n\in \Z$. For
each $V$ there exists a unique conformal mapping (the
quasi-momentum) $k:\cZ\to K(h) $ such that (see [Mi1])
$$
\cos k(z)= \D(z),\ \  z\in \cZ =\C\sm\cup \bar \g_n, \ \ \ \
K(h)=\C\sm\cup \G_n ,\ \  \  \G_n=(\pi n-ih_n,\pi n+ih_n),
$$
$$
\ \ \ k(z)=z+o(1)\qqq  as \qq |z|\to \iy ,
$$
where $\G_n$ is the vertical slit and the height $h_n\geq 0$ is
defined by the equation $\cosh h_n = (-1)^n\D(z_n) \geq 1$. We
emphasize that the introduction of the quasi-momentum $k(z)$
provides a natural labeling of all gaps $\g_n$ (including the empty
ones!) by demanding that $k(\cdot )$ maps the gap $g_0$ on the
vertical slit $\G_0=(-ih_0,ih_0)$.

For any $p\ge 1$ and  the weight $\o=(\o_n)_{n\in \Z}$, where
$\o_n\ge 1$, we introduce the real  spaces
$$
 \ell^p_{\o}=\{f=(f_n)_{n\in \Z}:\ \ \| f \|_{p,\o}<\iy \},\qq
\ \ \ \| f \|_{p, \o }^p=\sum_{n\in\Z } \o_n f_n^p <\iy .
$$
If the weight $\o_n=1$ for all $n\in \Z$, then we will write
$\ell^p_0=\ell^p$  with the norm  $\| \cdot \|_{p}$. For each $V$
 we introduce the sequences
$$
h=(h_n)_{n\in \Z},\qq  \g=(|\g_n|)_{n\in \Z}, \qq
 J=(J_n)_{n\in \Z},\qq J_n=|A_n|^{1\/2}\ge 0,\qq
A_n={2\/\pi }\int _{\g_n}v(z)dz\ge 0.
$$
For the defocussing cubic non-linear Schr\"odinger equation (a
completely integrable infinite dimensional Hamiltonian system),
$A_n$ is an action variables (see \cite{FM}). Recall the following
identities from \cite{KK1}
\[
\lb{id} {2\/\pi}\int_\R v(z)dz={1\/\pi}\iint_{\C
}|z'(k)-1|^2\,dudv={1\/2}\|V\|^2,
\]
where $\|V\|^2=\int_0^1(V_1^2(t)+V_2^2(t))dt$. Korotyaev \cite{K1}
obtained the two-sided estimates for the case $\ell^2$, for example
\[
\lb{ek1} \|g\|_2\le 2\|h\|_2\le\pi \|g\|_2(2+\|g\|_2^2),
\]
\[
\lb{ek2}
 {1\/\sqrt 2}\|g\|_2\le \|V\|\le 2\|g\|_2(1+\|g\|_2).
\]
 Our main goal is to obtain similar estimates in terms of the
$\ell^p$ and $\ell_\o^p-$norms.  Our first results are devoted to
estimates in terms of $\ell^p-$ norms. Let below ${1\/p}+{1\/q}=1,
p,q\ge 1$.

\begin{theorem}\lb{T1}
Let  $V\in L^1(0,1)$. Then the following estimates hold true:
\[
\lb{T1-1} 2^{-p}\|g\|_p\le  \|h\|_p\le 2\|g\|_p(1+\a_p^0\
\|g\|_p^p), \qqq \qqq p\in[1,2], \ \ \ \a_p^0={2^{(p+3)p}\/\pi} ,
\]
\[
\lb{T1-2} \|h\|_{p}\le \frac{2}{\pi}C_p^2 \|g\|_q\lt(
1+\lt({2C_p\/\pi^2}\rt)^{2\/p-1}\|g\|_q^{2\/p-1}\rt), \ \
C_p=\lt({\pi ^2\/2}\rt)^{1/p},  \ \  p\ge 2, \ \ {1\/p}+{1\/q}=1,
\]
\[
\lb{T1-3} {\|g\|_p  \/  2}\le  \|J\|_p \le {2\/  \sqrt{\pi
}}\|g\|_p(1+\a_p^0\|g\|_p^p)^{1/2},\qqq \ \ p\ge 1,
\]
\[
\lb{T1-4} {\sqrt{\pi}  \/  2}\|J\|_p\le \|h\|_p\le 4\|J\|_p(1+\a_p^0
2^{p}\|J\|_p^p),\qqq \ p\ge 1.
\]
\end{theorem}

Estimates \er{2.2}-\er{2.5} are new for $p\in [1,2)$.

In the case $z_n^-<z_n^+$ we define  the effective masses
$\mu_n^{\pm}$ by
\[
\lb{2.20} z(k)-z_n^{\pm} ={(k-\pi n)^2\/ 2\mu_n^{\pm}}(1+O(k-\pi n))
\qqq {\rm as} \qq \ z\to z_n^{\pm}.
\]
If $|g_n|=0$, then we set $\m_n^{\pm}=0$. Define the sequence
$\m^{\pm}=(\m_n^{\pm})_{n\in \Z}$. Our second result is

\begin{theorem}
\lb{T2} Let $h\in \ell_{\o}^{p}, p\in [1,2]$. Then the following
estimates hold true
\[
\lb{T2-1} \|h\|_{\iy }\le \min \rt\{2\pi \|\m^{\pm}\|_{\iy },\
\|J\|_{p,\o} ,\ 2\|g\|_{p,\o}(1+\a_p^0\|g\|_{p,\o}^p)^{1\/q}\lt\},
\]
\[
\lb{T2-2} \|g\|_{p,\o}\le 2\|h\|_{p,\o} \le c_0^9 \|g\|_{p,\o},\ \qq
c_0=e^{{1\/\pi}\|h\|_{\iy }},
\]
\[
\lb{T2-3} \|g\|_{p,\o} \le  2\|J\|_{p,\o} \le c_0^5 2\|g\|_{p,\o},
\]
\[
\lb{T2-4} {\sqrt{\pi} \/  2}\|J\|_{p,\o} \le \|h\|_{p,\o}\le
c_0^5\sqrt{\pi \/  2}\|J\|_{p,\o} ,
\]
\[
\lb{T2-5} \|g\|_{p,\o}\le 2\|\m^{\pm}\|_{p,\o} \le c_0^{18}
\|g\|_{p,\o}.
\]
\end{theorem}

Estimates \er{2.6}-\er{2.10} are new. Introduce  the real Hilbert
spaces $\ell_{(m)}^2, m\in \R$ of the sequences $\{f_n \}_1^{\iy }$
equipped with the norm $\| f \|_{(m)}^2=\sum _{n\ge 1}(2\pi n)^{2m
}f_n^2$. Korotyaev obtained the two-sided estimates for the
$\ell^2_{(m)}-$ norms, $m\ge 0$ for the even case $h_{-n}=h_n, n\in
\Z$ (\cite{K2}-\cite{K4}) and for $\ell^2_1-$norms without symmetry
(\cite{K1}, (\cite{K6})) (in all these estimates the factor
$c_0=e^{\|h\|_{\iy }/\pi}$ is absent).

\begin{proposition}
\lb{T4} Let $V\in L^2(0,1)$ . Then
\[
\lb{T4-1} \|h\|_\iy\le \|V\|,
\]
\[
\lb{T4-2}  \|V\|^2\le{2\/\pi}\|h\|_{p}^{}\,\|g\|_q, \ \ \  p\ge 1,
\]
\[
\lb{T4-3}
 \|V\|^2\le ({2\/\pi })^{2\/p}\|h\|_{p}^{2\/q}\,\|g\|_p^{2\/p},
 \ \ \qqq  p\in [1,2],
\]
\[
\lb{T4-4} \|V\|^2\le{2\/\pi}\|h\|_{\iy}^{}\,\|g\|_1\le
{4\/\pi^2}\|g\|^2_1,
\]
\[
\lb{T4-5} \|h\|_{\iy}^{}\le {2\/\pi}\|g\|_1,\qqq \|g\|_1\le2\|h\|_1.
\]
\end{proposition}

Note that the comb mappings are used in various fields of
mathematics. We  enumerate the more important directions:

{\it \no 1) the conformal mapping theory, 2) the L\"owner equation
and the quadratic differentials, 3) the electrostatic problems on
the plane, 4) analytic capacity, 5)  the spectral  theory of the
operators  with  periodic coefficients, 6) inverse problems for the
Hill operator and the Dirac operator, 7) the KDV equation and the
NLS equation with periodic initial value problem.}

{\bf Example of an electrostatic field.} Consider the system of
neutral conductors $\G_n, n\in \Z$  on the plane for some
$h\in\ell^p_\o$. In other words, we embed the system of neutral
conductors $\G_n$ in the external homogeneous electrostatic field
$E_0=(0, -1)\in \R^2$ on the plane. Then on each conductor there
exists the induced charge, positive $e_n>0$ on the lower half of the
conductor $\G_n$ and negative $(-e_n)<0$ on the upper half of the
conductor $\G_n$, since their sum equals zero. As a result we have
new perturbed electrostatic field $\cE\in \R^2$. It is well known
that $\cE=\ol {iz'(k)}=-\nabla y(k), \ k=u+iv\in K(h), \ \ z=x+iy$,
where $z(k)$ is the conformal mapping from $K(h)$ onto the domain
$\cZ=\C\sm \cup g_n$. The function $y(k)$ is called the potential of
the electrostatic field in $K(h)$. The density  of the charge on the
conductor has the form $\r(k)=|y_u'(k)|/4\pi, \ k\in \G_n$ (see
\cite{LS}). Thus we obtain the induced charge $e_n$ on the upper
half of the conductor $\G_n^+=\G_n\cap\C_+$ by:
$$
e_n={1 \/  4\pi}\int _{\G_n^+}x_v'(k)dv={1 \/  4\pi}|g_n|.
$$
Introduce the bipolar moment $d_n$ of the conductor $\G_n$ with the
charge density $\r(k)$ by $d_n={1\/4\pi}\int _{\G_n}vx_v(k)dv\ge 0$.
We transform this value into the form
$$
d_n={1 \/  2\pi}\int _{g_n}v(x)dx={A_n\/ 4}.
$$
In the paper \cite{KK3} we study inverse problems for both the
charge mapping $h\to e=(e_n)_{n\in \Z}$ and the bipolar moment
mapping  $h\to J$ acting in $\ell^p_{\o }, p\in [1,2]$. In order to
solve the inverse problems we need a priori estimates from Theorems
\ref{T1} and \ref{T2}.

A priori estimates of potentials in terms of spectral data
 essentially simplify the proof in the inverse
problems. Such simplification was introduced by Garnett and
Trubowitz \cite{GT1} and Kargaev and Korotyaev \cite{KK2} and
essentially was used in \cite{K6}-\cite{K9}.

For the sake of the reader, we briefly recall the results existing
in the literature about the a priori estimates. Firstly, we describe
a priori estimates for the Hill operators $-{d^2\/dt^2}+P$ in
$L^2(\R)$ with the real 1-periodic potential $P\in L^2(0,1)$ nd
$\int_0^1P(t)dt=0$. The spectrum of this operator consists of
intervals separated by gaps $\g_n, n\ge 1$ with the lengths
$|\g_n|\ge 0, n\in \Z$. Marchenko and Ostrovski [MO1-2] obtained the
estimates: $\|P\|\le C(1+\|h\|_{\iy})\|h\|_{(1)},\ \ \|h\|_{(1)}\le
C\|P\|e^{C\|P\|}$ for some absolute constant $C$, where
$\|P\|^2=\int_0^1P^2(t)dt$ and $\|h\|_{(1)}^2=\sum (2\pi n)^2h_n^2$.
These estimates are very rough since they used the Bernstein
inequality. Using the harmonic measure argument Garnett and
Trubowitz \cite{GT} obtained $\|\g\|\le (4+\|h\|_{(1)})\|h\|_{(1)}$,
where $\g=(|\g_n|)_1^\iy$ and recall that $h_n$ are heights  on the
quasimomentum domain. First two-sided estimates (very rough) for $g,
h$ were obtained in \cite{KK2}.

Identities and a complete system of a priori estimates  (in terms of
gap lengths, effective masses etc) were obtained by Korotyaev
\cite{K1}-\cite{K6}. In the paper [K3], [K5], [K6] the following
estimates  were obtained:
$$
\|\g\|_2\le 6\|P\|(1+\|P\|^{1\/3}) ,\qqq \|P\|\le
4\|\g\|_2(1+\|\g\|_2^{1\/3}) ,
$$$$
2\|h\|_{(1)}\le \pi \|P\|( 1+\|P\|^{1\/3}) ,\ \ \ \ \|P\|\le
3(6+\|h\|_\iy)^{1\/2}\|h\|_{(1)}.
$$
These estimates show the ''equivalence'' of the values $\|\g
\|,\|h\|_{(1)}, \|P\|$. Note the author solved the inverse problem
and obtained the two-sided estimates for the case $P=y'$, where
$y\in L^2(0,1)$. A priori two-sided estimates for the case
$P^{(m)}\in L^2(0,1), m\ge 1$ were obtained by Korotyaev  in
\cite{K1}.

Secondly for the Zakharov-Shabat systems the two-sided estimates
were obtained by Korotyaev for $V\in L^2(0,1)$ in \cite{K1} (see
\er{ek1}-\er{ek2}) and for $V'\in L^2(0,1)$ in \cite{K6}. There are
no a priori estimates for the case $V^{(m)}\in L^2(0,1), m\ge 2$.
 The proof of a priori estimates in \cite{K1}-\cite{K6} is based on
 the analysis of the quasimomentum as the conformal mapping.

We shortly describe the proof. In order to prove Theorem \ref{T1}
-Proposition \ref{T4} we use the analysis of a conformal mapping
corresponding to quasimomentum of the Zakharov-Shabat operator. That
makes it possible to reformulate the problems for the differential
operator as the problems of the conformal mapping theory. Then we
should study the metric properties of a conformal mapping from
$\C_+$ onto a "comb" $K_+(h)$. A similar analysis was done partially
in [KK1-3], [K1-3]. In the present paper we use an approach, based
on the identities for the Dirichlet integral (1.2) from \cite{KK1}
and the estimates from Theorem 2.7 and Lemma 3.1. We emphasize the
important role of the Dirichlet integral (this is the energy for the
conformal mapping) in this consideration.

We now describe the plan of the paper. In Section 2 we shall obtain
some preliminaries results and "local basic estimates" in Theorem
\ref{T3.5}. In Section 3 we shall prove Lemma 3.1 and the main
theorems. Moreover, we consider some examples, which describe our
estimates.

\section{Preliminaries}
\setcounter{equation}{0}

Consider a conformal mapping $z:K_+(h)\to\C_+$ with asymptotics
$z(iv)=iv(1+o(1))$ as $v\to\iy$, where  $k=u+iv\in K(h)$. Here
$h=(h_n)_{n\in\Z}\in\ell^{\iy}, h_n\ge 0$ is some sequence and
$K_+(h)=\C_+\cap K(h)$ is the so-called comb domain, where $K(h)$ is
given by
$$
K(h)=\C \sm \cup_{n\in\Z}\G_n,\qq \G_n=[u_n -ih_n,u_n+ih_n], \qq
u_*=\inf_n (u_{n+1}-u_n)\ge 0,
$$
where $u_n, n\in \Z$ is strongly increasing sequence of real numbers
such that $u_n\to \pm \iy $ as $n\to \pm \iy $. We fix the sequence
$u_n, n\in \Z$ and consider the conformal mapping for various $h\in
\ell^\iy$. For fixed $h$ the difference of any two such mappings
equals a real constant, but the imaginary part $y(k)=\Im z(k)$ is
unique. We call such mapping $z(k) $ the comb mapping. Define the
inverse mapping $k(\cdot): {\C _+}\to K_+(h)$. It is clear that
$k(z), z=x+iy\in \C _+$ has the continuous extension into $\ol \C
_+$. We define "gaps" $g_n$, "bands" $\s_n$ and the "spectrum" $\s$
of the comb mapping by:
$$
g_n=(z_n^-, z_n^+)=(z(u_n-0), z(u_n+0)),\qqq \s_n=[z_{n-1}^+,
z_n^-],\qqq \s=\cup_{n\in\Z }\s_n.
$$
The function $u(z)=\Re k(z) $ is strongly  increasing on each band
$\s_n$ and $u(z)=u_n$  for all $z\in [z_n^-,z_n^+],\ n\in\Z $; the
function $v(z)=\Im k(z)$ equals zero  on each  band $\s_n$ and is
strongly convex on each gap $g_n\ne \es$ and has the maximum  at
some point $z_n $ given by $v(z_n)=h_n$. If  the gap  is empty we
set $z_n=z_n^{\pm}$. The function $z(\cdot) $ has an analytic
extension (by the symmetry) from  the domain $K_+(h)$   onto  the
domain $K(h)$ and $z(\cdot): K(h)\to z(K(h))=\cZ=\C\sm\cup \ol g_n $
is  a conformal mapping. These and  others properties of the comb
mappings it is possible  to  find in the papers of  Levin \cite{Le}.

We formulate our second result about the estimates for conformal
mappings.

\begin{theorem}\lb{T2.1}
Let  $u_*=\inf_n (u_{n+1}-u_n)>0$. Then the following estimates hold
true
\[
\lb{2.2} \|h\|_p\le 2\|g\|_p(1+\a_p\ \|g\|_p^p), \qqq \qq p\in
[1,2], \ \ \ \a_p={(2+\pi )^p2^{p(p+2)}\/\pi u_*^p} ,
\]
\[
\lb{2.3} \|h\|_{p}\le {2\/\pi}C_p^2 \|g\|_q\lt( 1+\lt[{2C_p\/\pi
u_*}\rt]^{{2\/p-1}}\|g\|_q^{2\/p-1}\rt), \ \ C_p=\lt({\pi
^2\/2}\rt)^{1/p},  \ \  p\ge 2, \ \ {1\/p}+{1\/q}=1,
\]
\[
\lb{2.4} {\|g\|_p  \/  2}\le  \|J\|_p \le {2\/  \sqrt{\pi
}}\|g\|_p(1+\a_p\|g\|_p^p)^{1/2},\ \ \ \ \ p\ge 1,
\]
\[
\lb{2.5} {\sqrt{\pi}  \/  2}\|J\|_p\le \|h\|_p\le 4\|J\|_p(1+\a_p
2^{p}\|J\|_p^p),\qqq \ \ \ p\ge 1.
\]
\end{theorem}

We formulate our second result about the estimates for conformal
mappings.

\begin{theorem}
\lb{T2.2} Let $h\in \ell_{\o}^{p}, p\in [1,2]$ and let $u_*>0$. Then
the following estimates hold true
\[
\lb{2.6} \|h\|_{\iy }\le \min \{2\pi \|\m^{\pm}\|_{\iy },\
\|J\|_{p,\o} ,\ 2\pi
^{-1\/p}\|g\|_{p,\o}(1+\a_p\|g\|_{p,\o}^p)^{1/q}\},
\]
\[
\lb{2.7} \|g\|_{p,\o}\le 2\|h\|_{p,\o} \le c^9 \|g\|_{p,\o},\ \qq
c=e^{\|h\|_{\iy }/u_*},
\]
\[
\lb{2.8} \|g\|_{p,\o} \le  2\|J\|_{p,\o} \le c^5 2\|g\|_{p,\o},
\]
\[
\lb{2.9} {\sqrt{\pi} \/  2}\|J\|_{p,\o} \le \|h\|_{p,\o}\le
c^5\sqrt{\pi \/  2}\|J\|_{p,\o} ,
\]
\[
\lb{2.10} \|g\|_{p,\o}\le 2\|\m^{\pm}\|_{p,\o} \le c^{18}
\|g\|_{p,\o}.
\]
\end{theorem}

Define the Dirichlet integral $I_D$ and the moment $Q_0$ by
$$
I_D={1\/\pi}\iint_{\C }|z'(k)-1|^2\,dudv ={1\/\pi}\iint_{\C
}|k'(z)-1|^2\,dxdy,\qqq \qq Q_0={1\/\pi}\int_\R v(z)dz,
$$
where $k=u+iv,\ z=x+iy$. The last identity holds since the Dirichlet
integral is invariant under the conformal mappings.  In order to
prove our main theorems we need the following

\begin{proposition}
\lb{T3.6} Let $h\in \ell^{\iy}$ and let $u_*\ge 0$. Then
\[
\lb{2.29} {\|h\|_\iy^2\/2}\le Q_0,
\]
\[
\lb{3.17} \pi Q_0\le\|h\|_{p}^{}\,\|g\|_q, \ \ \  p\ge 1,
\]
\[
\lb{3.18}
 I_D\le ({2\/\pi })^{2\/p}\|h\|_{p}^{2/q}\,\|g\|_p^{2/p},
 \ \ \qqq  p\in[1,2],
\]
\[
\lb{3.19} \pi Q_0\le\|h\|_{\iy}^{}\,\|g\|_1\le {2\/\pi}\|g\|^2_1,
\]
\[
\lb{3.20} \|h\|_{\iy}^{}\le {2\/\pi}\|g\|_1,\qqq \|g\|_1\le2\|h\|_1.
\]
\end{proposition}

{\bf Proof of Theorem \ref{T1}- Proposition \ref{T4}} follow
directly from Theorem \ref{T2.1}- Proposition \ref{T3.6} and the
identity \er{id}.

We recall needed results.
Below we will use very often the following simple estimate
\[
\lb{1.3} |g_n|\le 2h_n,\qq {\rm all}\qq n\in\Z,
\]
see e.g. \cite{MO1}, \cite{KK1}. Hence if $h\in \ell^p_{\o}$, then
$\g\in \ell^p_{\o} $.

For each $h\in \ell^{\iy}$ the following estimates and identities
hold true
\[
\lb{1.5} {1 \/4} \|g\|^2\le 2Q_0=I_D=\sum A_n=\|J\|^2\le {2\/\pi}
\sum_{n\in\Z }h_n|g_n|,
\]
\[
\lb{2.30} \max \lt\{ {|g_n|^2 \/  4}, {|g_n|h_n \/  \pi}\rt\}\le
A_n= {2 \/  \pi} \int _{g_n}v(x)dx\le {2|g_n|h_n \/  \pi},
\]
see [KK1]. These show that functional $Q_0={1\/\pi}\int_{\R }v(x)\,dx $
is bounded for $h\in \ell^2$.

Below we  will sometimes write $g_n(h), z(k,h),..$, instead of $g_n,
z(k),..$, when several sequences $h\in \ell^\iy$ are being dealt
with. Recall the Lindel\"of principle (see [J]), which is formulated
in  the form, convenient  for us (see \cite{KK1}):

{\it  Let $h,\ \wt h\in \ell^\iy;$ and let $\wt h_n\le h_n$ for all $n\in\Z $.
Then the following estimates hold:}
\[
\lb{2.23} y(k,\wt h)\ge y(k,h),\ \qqq all \qq k\in K_+(h),
\]
\[
\lb{2.26} Q_0(\wt h)\le Q_0(h)\ \qqq  and \ if \qq Q_0(\wt h)=
Q_0(h),\ \ then \ \wt h=h,
\]
\[
\lb{2.27} |\s_n(\wt h)|\ge |\s_n(h)|.
\]

Define the effective masses $\n_n$ in the plane $K(h)$ for the  end
of the slit $[u_n+ih_n,u_n-ih_n], h_n>0$ by
\[
\lb{2.21} k(z)-(u_n+ih_n)={(z-z_n)^2 \/  2i\n_n}(1+O(z-z_n))) \qq as
\qq z\to z_n.
\]
Thus we obtain $\n_n=1/|k''(z_n)|$, if $h_n>0$ and we set $\n_n=0$
if $|g_n|=0$.  We show the possibility of the Lindel\"of principle
in the following Lemma.

\begin{lemma} \lb{L2.5}
For each $h\in\ell^\iy$ the estimate \er{2.29} and the following
estimate hold true
\[
\lb{2.28}
\nu_n\le h_n, \qqq all \qq n\in\Z.
\]
\end{lemma}
\no {\it Proof.} Sufficiently to proof for the case $n=0$. We apply
estimate \er{2.23} to $h$ and to the new sequence:
 $\wt h_0=h_0$ and $\wt h_n=0$ if $ n\ne 0$.  It  is clear
that  $z(k,\wt h)=\sqrt{(k-u_0)^2+h_0^2}$ (the principal value).
Then \er{2.23} gives
$$
y(k,h)\le\Im(\sqrt{(k-u_n)^2+h_0^2}),\ k\in K_+(h).
$$
Then asymptotics \er{2.21} of the function $z(k,h)$ as
 $k\to u_0+ih_0$ yields \er{2.28}.
In  order  to  prove \er{2.29} we use \er{2.26} since $Q_0(\wt
h)=h_0^2/2$. Note that it and  \er{2.26} yield \er{2.29}.
 $\BBox$

We recall estimates from [K4].

\begin{theorem} \lb{T3.3}
 Let $h\in \ell^{\iy}$. Then for any  $r>0,n\in\Z$ the following estimate holds
 true
\[
\lb{3.3} h_n^2\le{\pi\/4}
\max\lt\{1,{h_n\/r}\rt\}\iint_{u_n+S_r}|z'(k)-1|^2dudv,\qq
S_r=\{z\in\C :\ |\Re z|<r\},
\]
If in addition, $\inf_n ( u_{n+1}-u_n)=u_*>0$, then
\[
\lb{3.6} {\pi\/4} I_D\le\|h\|_2^2\le {\pi^2\/2}\max
\rt\{1,{\|h\|_\iy\/u_*}\rt\} I_D\le {\pi^2\/2}\max
\rt\{1,{I_D^{1/2}\/u_*}\rt\}I_D,
\]
\[
\lb{3.7} {1\/2}\|g\|_2\le\|h\|_2\le\pi
\|g\|_2\rt(1+{2\/u_*^2}\|g\|_2^2\rt),
\]
\[
\lb{3.8} {\|g\|_2\/2}\le \|J\|_2\le
\sqrt{2}\|g\|_2(1+{\sqrt{2}\/u_*}\|g\|).
\]
\end{theorem}

 In order to prove Theorem \ref{T3.5} we need the following result about
the simple mapping and the domain $ S_r=\{z\in\C :\ |\Re z|<r\}, \
r>0$.

\begin{lemma} \lb{L3.4}
The function $f(k)=\sqrt{k^2+h^2},\ k\in\C \sm [-ih,ih],
h>0\,$ is the conformal mapping from $\C\sm [-ih, ih]$ onto
$\C\sm [-h,h]$ and $S_r\sm [-h,h]\subset f(S_r\sm [-ih,ih])$  for any $r>0$.
\end{lemma}
\no {\it Proof.} Consider the image of the half-line $k=r+iv,
v>0$. We have the equations
\[
\lb{3.9}
x^2+y^2=\x \ev r^2+h^2-v^2,\ \ \ \ \  \ \ xy= rv.
\]
The second identity in \er{3.9} yields $x>0$ since $y>0$ . Then
$x^4-\x x^2-r^2v^2=0,$
and enough to check the following inequality
$x^2={1\/2}(\x +\sqrt{\x ^2+4r^2v^2})  > r^2.$
The last estimate follows from the simple relations
$$
(r^2+h^2-v^2)^2+4r^2v^2>(r^2+v^2-h^2)^2,\ \ \ \ \ \ 4r^2v^2>4r^2(v^2-h^2).
\ \ \ \ \BBox
$$

We prove the local estimates for the small slits, which are crucial
for us.

\begin{theorem} \lb{T3.5}
Let $h\in \ell^{\iy}$. Assume that
$(u_n-r, u_n+r)\subset(u_{n-1}, u_{n+1})$ and $h_n\le {r\/2}$,
for some $n\in\Z $ and $r>0$. Then
\[
\lb{3.10}
|h_n-|\mu_n^{\pm}||\le \frac{2+\pi}{r}|\mu_n^{\pm}|\sqrt{I_n},
\ \ \ I_n=\frac1{\pi}\int\!\!\!\!\int_{u_n+S_{r}}|z'(k)-1|^2dudv,
\]
\[
\lb{3.11}
0\le h_n-\nu_n\le 2\frac{2+\pi}{r}h_n \sqrt{I_n},
\]
\[
\lb{3.12} 0\le h_n-\frac{|g_n|}2\le \frac{2+\pi}{r}h_n\sqrt{I_n}.
\]
\end{theorem}
\no {\it Proof.} Define the functions $f(k)=\sqrt{k^2+h_n^2}$, $\
k\in S_{r}\sm [-ih_n,ih_n], \f=f^{-1}$
 and $F(w)=z(u_n+\f(w),h), w=p+iq$, where the variable
$w\in G_1=f\big((S_{r}\sm [-ih_n,ih_n]\big))$. The function $F$ is
real for real $w$, then $F$ is analytic in the domain $G=G_1\cup
[-h_n,h_n]$ and Lemma \ref{L3.4} yields $S_{r}\ss G$. Let now
$|w_1|={r\/2}$ and $B_r=\{z:|z|<r\}$. Then the following estimates
hold
\[
\lb{3.13}
 \sqrt{\pi}\frac{r}2|F'(w_1)-1| \le (\iint_{B_r}
\big|F'(w)-1\big|^2dpdq)^{1/2}\le
\]
$$
 \le\bigg(\iint_{B_r}|(F(w)-\f(w))'|^2dpdq\rt)^{1/2}+\rt(\iint_{B_r}
\big|\f'(w)-1\big|^2dpdq\rt)^{1/2}.
$$
The invariance of the Dirichlet integral with respect to
the conformal mapping gives
\[
\lb{3.14}
 \iint_{B_r}\!\!|(F(w)-\f(w))'|^2dpdq=
\iint_{\f(B_r)}\!\!\!\!|z'(k)-1|^2dudv\le
\iint_{S_{r}+u_n}\!\!|z'(k)-1|^2dudv=\pi I_n.
\]
Moreover, the identity $2Q_0=I_D$ implies
\[
\lb{3.15} {1\/\pi}\iint_{B_r}|\f'(w)-1|^2dpdq\le
\frac1{\pi}\iint_{\C }|\f'(w)-1|^2dpdq
=\frac2{\pi}\int_{-h_n}^{h_n}\sqrt{h_n^2-x^2}\ dx=h_n^2.
\]
Then \er{3.13}-\er{3.15} for $|w_1|={r\/2}$  yields
$
\big|F'(w_1)-1\big|\le\frac2{r}(\sqrt{I_n}+h_n),
$
and \er{3.3} gives
$$
h_n^2\le {\pi^2\/4}\cdot{1\/\pi}\iint_{S_{r}+u_n}|z'(k)-1|^2dudv\le
\frac{\pi^2}4 I_n.
$$
Then for $|w_1|={r\/2}$  we have
\[
\lb{3.16}
|F'(w_1)-1|\le\frac2{r}(1+\frac{\pi}2)\sqrt{I_n}={2+\pi \/  r}\sqrt{I_n},
\]
and the maximum principle yields the needed estimates for $|w_1|\le r/2.$

We prove \er{3.10} for $\mu_n^{+}$. The definition of $\m_n^{\pm}$
(see \er{2.20}) implies
$$
F'(h_n)=\lim_{x\searrow h_n}z'(u_n+g(x))\cdot
g'(x)=\lim_{x\searrow h_n}\frac{g(x)}{\mu_n^{+}}\cdot
\frac x{g(x)}=\frac{h_n}{\mu_n^{+}}.
$$
The substitution of the last identity into \er{3.16} gives \er{3.10}.
The proof for $\m_n^-$ is similar.

We show \er{3.11}.  The definition of  $\nu_n$ (see \er{2.21}) yields
$$
(z(k)-z_n)^2=2i\nu_n(k-u_n-ih_n)(1+o(1)) \ \ as \  \ k\to u_n+ih_n,
$$
$$
\f(w)-ih_n=-\frac i{2h_n}(w-z_n)^2(1+o(1))\ \ as \ w\to u_n.
$$
Then we have $F'(0)=\sqrt{\frac{\nu_n}{h_n}}$ and the substitution
of the last identity into \er{3.16} shows
$\bigg|\sqrt{\frac{\nu_n}{h_n}}-1\bigg|\le {2+\pi\/r} \sqrt{I_n}$,
which gives \er{3.11}, since by \er{2.30}, $\nu_n\le h_n$. Estimate
\er{3.16} yields
$$
0\le2
h_n-|g_n|=\int_{-h_n}^{h_n}(1-F'(x))dx\le\frac{2h_n}r(2+\pi)\sqrt{I_n},
\qqq
$$
which implies \er{3.12}.
\BBox

 We prove the estimates in terms of the $\ell^p$-norms.

\no {\bf  Proof of Proposition  \ref{T3.6}.} We have proved
\er{2.29} in Lemma \ref{L2.5}. Estimate \er{1.5} and the H\"older
inequality yield \er{3.17}. Using \er{2.29},  \er{1.5} and the
H\"older inequality, we obtain
$$
(\pi /2) I_D=\pi Q_0\le\sum h_n|g_n|\le\|h\|_{\iy}^{1-{p \/
q}}\sum_{n\in\Z }|g_n| \,h_n^{p/q}\le(I_D)^{(1-{p \/
q})/2}\|g\|_p\|h\|_p^{{p \/ q}},
$$
$$
 (\pi /2) I_D^{\frac{(1+{p\/q})}2}\le \|g\|_p\|h\|_p^{{p\/q}},
\ \ \ \ \ {\rm and}\ \ \ \ \ \ I_D\le
\rt({2\/\pi}\rt)^{{2\/p}}\,\|g\|_p^{{2\/p}} \|h\|_p^{{2\/q}}.
$$
Estimate \er{3.17} at $q=1$ implies the first one in \er{3.19}. The
last result and \er{2.29} yield the first inequality in \er{3.20}
and then the second one in \er{3.19}. The second estimate in
\er{3.20} follows from $|g_n|\le2 h_n, \ n\in\Z $ (see \er{1.3}).
 $\BBox$

\section{Proof of the mains theorems}
\setcounter{equation}{0}

\no {\bf Proof of Theorem \ref{T2.1}}. Let $p\in [1,2]$ and
$r={u_*\/2}$. Estimate \er{3.3} implies
\[
\lb{3.21} h_n^2\le{\pi^2\/4} \max\{1,{h_n \/  r}\}I_n,  \ \ \
I_n={1\/\pi} \iint_{u_n+S_r}\!\!\!|z'(k)-1|^2dudv,\ \ \
\]
Hence
\[
\lb{3.22}
  h_n\le {\pi^2\/u_*}I_n,\ \ \ if \ \  h_n>{u_*\/4},
\ \ and \ \
  h_n\le {\pi \/2}\sqrt{I_n},\ \ \ if \ \  h_n\le {u_* \/4}.
\]
Moreover, \er{3.12} yields
\[
\lb{3.23} h_n\le {|g_n|\/2}+2 {2+\pi\/u_*}h_n\sqrt{I_n},\qqq if \qq
h_n<{u_*\/4},
\]
and then
\[
\lb{3.24} h_n\le 2\pi \frac{2+\pi }{u_*}I_n, \qqq if\qq
h_n<{u_*\/4}, \ |g_n|\le h_n,
\]
since $h_n\le {\pi \/2}\sqrt{I_n}$. Hence using \er{3.24},
\er{3.22}, we obtain
$$
if\ \  \ |g_n|\le h_n \ \  \Rightarrow \ \
  h_n\le 2\pi \frac{2+\pi}{u_*}I_n=C_1I_n,\ \ \ \ C_1=2\pi {2+\pi\/u_*}.
$$
The last inequality and \er{2.29} yield
\[
\lb{3.25} \| h\|_p\le (\sum _{h_n<|g_n|}h_n^{p})^{1\/p}+ (\sum
_{|g_n|\le h_n}h_nh_+^{p-1})^{1\/p} \le \|g\|_p+
C_1^{\frac1{p}}I_D^{\frac{p+1}{2p}},\qq h_+=\|h\|_\iy.
\]
If we assume that $C_1^{{1\/p}}I_D^{{p+1\/2p}}\le \|g\|_p$, then we
obtain $\|h\|_p\le 2\|g\|_p.$

Conversely,  if we assume that $\|g\|_p\le
C_1^{\frac1{p}}I_D^{\frac{p+1}{2p}}$, then \er{3.25}, \er{3.18}
implies
$$
\|h\|_p\le 2C_1^{{1\/p}} \lt[\rt({2\/\pi} \rt)^{2\/p}
\|h\|_{p}^{2\/q}\,\|g\|_p^{2/p}\rt]^{{p+1\/2p}}.
$$
Hence
$$
\|h\|_p^{1/p^2}\le 2C_1^{1\/  p}\rt[(2 /  \pi)\|g\|_p\rt]^{p+1 \/
p^2} \ \ \ \ \ \Rightarrow \ \ \ \ \ \ \|h\|_p\le 2^{p^2}C_1^p
(2/\pi) ^{p+1} \|g\|_p^{1+p},
$$
which yields \er{2.2}.

Let $p\ge 2$. Using inequality \er{3.6}, \er{2.29} we obtain
\[
\lb{3.26} \|h\|_p\le (\sum h_+^{p-2} h_n^2)^{1\/p}\le
C_pb^{1\/p}I_D^{1\/2}, \ \ \ \ b=b(h_+),\ b(t)=\max\{1,{t\/r}\}, \
h_+=\|h\|_\iy.
\]
Consider the case $b\le 1.$ Then \er{3.26}, \er{3.17} imply
$$
\|h\|_p^2\le C_p^2I_D\le C_p^2(2/\pi) \|h\|_p\|g\|_q, \ \ \ \ \ {\rm
and} \ \ \ \ \ \ \ \|h\|_p\le (2C_p^2/\pi) \|g\|_q,
$$
Consider the case $b>1.$ Then the substitution of \er{2.29}, \er{3.17}
into \er{3.26} yield
$$
\|h\|_p\le C_pI_D^{\frac{p+1}{2p}}u_*^{-1/p}\le
u_*^{-1/p}C_p[(2/\pi) \|h\|_p\|g\|_q]^{\frac{p+1}{2p}},
$$
and
$$
\|h\|_p^{{p-1\/2p}}\le C_pu_*^{-1/p}[(2/\pi)
\|g\|_q]^{\frac{p+1}{2p}} \ \ \ \ \ \Rightarrow \ \ \  \ \
\|h\|_p\le (C_pu_*^{-1/p})^{\frac{2p}{p-1}}
(2/\pi)^{\frac{p+1}{p-1}}
 \|g\|_q^{\frac{p+1}{p-1}},
$$
and combining these two cases we have \er{2.3}.

 Estimate $|g_n|\le 2J_n$ (see \er{2.30}) yields the first one in \er{2.4}. Relation
\er{2.30} implies
$$
\|J\|_p^p=\sum |J_n|^p\le \sum (2/\pi )^{p/2}h_n^{p/2}
|g_n|^{p/2}\le (2/ \pi )^{p/2} \|h\|_p^{p/2} \|g\|_p^{p/2},
$$
and using \er{2.2} we obtain the second estimate in \er{2.4}:
$$
\|J\|_p\le \sqrt{2/\pi }\|g\|_p^{1/2}
[2\|g\|_p(1+\a_p\|g\|_p^p)]^{1/2}= {2\/  \sqrt{\pi
}}\|g\|_p(1+\a_p\|g\|_p^p)^{1/2},
$$
recall that  $\a_p=(2^{p+2}(2+\pi )/u_*)^p/\pi $. Inequality
$J_n^2\le 4h_n^2/\pi $ (see \er{2.30})    yields the first estimate
in \er{2.5}. Using \er{2.2} and $\|g\|_p\le 2\|J\|_p$ (see \er{2.4})
we deduce that
$$
\|h\|_p\le 2\|g\|_p(1+\a_p\|g\|_p^p)\le 4\|J\|_p(1+\a_p
2^{p}\|J\|_p^p). \ \ \ \BBox
$$

Recall the following identity for $v(z)=\Im k(z), z=x+iy$ from
[KK1]:
\[
\lb{3.32} v(x)=v_n(x)\big(1+Y_n(x)\big), \quad
Y_n(x)=\frac1{\pi}\!\!\! \intl_{\R \bs
g_n}\frac{v(t)dt}{|t-x|\,v_n(t)}, \quad
v_n(x)=|(x-z_n^{+})(x-z_n^{-})|^{1\/2},
\]
for all $x\in g_n=(z_n^{-},z_n^{+})$. In order to prove Theorem
\ref{T2.2} we need the following results.

\begin{lemma} \lb{L3.8}
Let $h\in \ell^\iy $ and $u_*>0$ and $c=e^{\|h\|_{\iy}\/u_*}$. Then
the following estimates hold:
\[
\lb{3.33} s=\inf |\s_n|\le u_* \le {\pi s\/2}\max\lt\{e^2,c^{5\pi
\/2}\lt\},
\]
\[
\lb{3.34} 1+{2\|h\|_{\iy}\/  s\pi}\le c^9 ,
\]
\[
\lb{3.35} \max_{n\in g_{n} } Y_n(x)\le {2\|h\|_{\iy}\/  \pi s}, \ \
n\in \Z,
\]
\[
\lb{3.36} 2h_n\le |g_n|(1+\max_{n\in g_{n} } Y_n(x))\le
|g_n|(1+{2\|h\|_{\iy} \/  \pi s}) \le |g_n| c^9, \ \ n\in \Z.
\]
\end{lemma}
\no {\it Proof.} Introduce the domain
$$
G=\big\{z\in\C : h_+\ge\Im\,z>0,\ \Re\,z
\in(-{u_* \/  2},{u_* \/  2})\big\}\cup\{\Im\,z>h_+\big\}, \ \
h_+=\|h\|_{\iy}.
$$
Let $F$ be the conformal mapping from $G$ onto $\C _+$, such that
$F(iy)\sim iy$ as $y\nearrow +\iy$ and let $\a , \b $ be images of
the points ${u_*\/2}, {u_*\/2}+ih_+$ respectively. Define the
function $f=\Im F$. Fix any $n\in \Z $. Then the maximum principle
yields
$$
y(k)=\hbox{Im}\,z(k,h)\ge f(k-p_n),\quad k\in G+p_n,\ \ \
p_n={1 \/  2}(u_{n-1}+u_n).
$$
Due to the fact that these positive functions equal zero
on the interval $(p_n- {u_*\/2},p_n+{u_*\/2})$, we obtain
$$
\pa_v y(x)=\pa_u x(x)\ge \pa_v f(x-p_n),\qqq x\in
\big(p_n-\frac{u_*}2, p_n+\frac{u_*}2\big).
$$
where $\pa_x={\pa \/\pa x}$. Then
$$
z(u_{n})-z(u_{n-1})\ge\intl_{-u_*/2}^{u_*/2}\pa_v f(x)\,dx= 2\a
>0,
$$
and the estimate $s\le u_*$ (see [KK1]) implies $2\a\le s\le u_*$.
Let $w: \C_+\to G$ be the inverse function for $F$, which is defined
uniquely and the Christoffel-Schwartz formula  yields
$$
w(z)=\int _0^z\sqrt{{t^2-\b^2 \/  t^2-\a^2}}dt,\ \ \ \ \ 0<\a<\b.
$$
Then we have
\[
\lb{3.37}
{u_* \/  2}=\int _0^{\a}\sqrt{{\b^2 -t^2\/  \a^2-t^2}}dt,
\ \ \
h_+=\int _{\a}^{\b}\sqrt{{\b^2-t^2 \/  t^2-\a^2}}dt.
\]
The first integral in \er{3.37} has the simple double-sided estimates
$$
\a=\int _{0}^{\a}dt\le {u_* \/  2}\le
\int _0^{\a}{\b dt\/  \sqrt{\a^2-t^2}}={\b\pi \/  2},
$$
that is
\[
\lb{3.38}
2\a\le s\le u_*\le \pi \b.
\]
Consider the second integral in \er{3.37}. Let $\ve =\b /\a\ge 5$ and using
the new variable $t=\a \cosh r, \cosh \d =\ve $, we obtain
$$
h_+=\a \int _0^{\d }\sqrt{\ve ^2-\cosh ^2r}dr\ge
\a\ve\int_0^{\d /2}\sqrt{1-{\cosh ^2r \/  \ve ^2}}dr\ge \b\d{2\/  5},
$$
since for $r\le \d/2$ we have the simple inequality
$$
{\cosh ^2r \/  \cosh ^2\d}  \le e^{-\d} (1+e^{-\d})^2 \le \ve^{-1}
(1+\ve^{-1})^2.
$$
Due to $\ve \le e^{\d}$ we get $\ve \le \exp (5h_+/2\b)$
and estimate \er{3.38} implies
\[
\lb{3.39}
{1 \/  s}\le {\pi  \/  2u_*}\exp ({5\pi  \/  2u_*}h_+),
\ \ \ if  \ \ \ \ve \ge 5.
\]
If $\ve \le 5$, then using \er{3.38} again we obtain
$$
{1 \/  s}\le {\ve  \/  2\b}\le {\pi \ve  \/  2u_*}\le
{\pi \/  2u_*}5, \ \ \ \ \ \    if  \ \ \ \ve \le 5.
$$
and the last estimate together with \er{3.39} yield \er{3.33},\er{3.34}.

 Identity \er{3.32} for $ x\in g_n=(z_n^{-},z_n^{+})$ implies
$$
\pi Y_n(x)= \int_{-\iy}^{z_n^--s}\frac{v(t)dt}{|t-x|v_n(t)} +
\int_{z_n^++s}^{\iy}\frac{v(t)dt}{|t-x|v_n(t)}\le
\int_{-\iy}^{z_n^--s}\frac{h_+dt}{|t-z_n^-|^2} +
\int_{z_n^++s}^{\iy}\frac{h_+dt}{|t-z_n^+|^2}\le
 {2h_+\/s}.
$$
Using \er{3.32}, \er{3.34} and simple inequality $v_n(z_n)\le
|g_n|/2$ we have \er{3.36}.  \ \ $\BBox$

We prove the two-sided estimates  of $h_n, |g_n|, \m_n^{\pm}, J_n$
in the weight spaces.

\no {\bf  Proof of Theorem \ref{T2.2}}.
 The first estimate in \er{2.6} follows from
$h_n\le 2\pi |\m_{n}^{\pm}|$ (see [KK1]). The second one in \er{2.6} follows from
 $\|h\|_{\iy}\le \sqrt{I_D}=\|J\|_2\le \|J\|_p\le \|J\|_{p,\o}$
since $\o_n\ge 1$ for any $n\in \Z$.
Moreover, substituting \er{3.18} into $\|h\|_{\iy}\le \sqrt{I_D}$,
using \er{2.2} and $\|f\|_p\le \|f\|_{p,\o}$  for any $f$,  we obtain  the last estimate in \er{2.6}.

Recall that $c=\exp {{\|h\|_{\iy }\/u_*}}$. The first estimate in
\er{2.7} follows from \er{1.3}. Due to \er{3.36} we  get $2h_n\le
c^9|g_n|$, which yields the second estimate in \er{2.7}.

The first estimate in \er{2.8} follows from \er{2.30}. Using
\er{2.30}, \er{3.36} we  have $J_n^2\le 2|g_n|h_n/\pi \le (c^9/\pi
)|g_n|^2$, which gives the second estimate in \er{2.8}.

The first estimate in \er{2.9} follows from \er{2.30}, \er{1.3}.
Using \er{2.30}, \er{3.36} we obtain $h_n^2\le c^9|g_n|h_n/2 \le
(\pi c^9/2)J_n^2$, which yields the second inequality in  \er{2.9}.

Identity $2|\m_n^{\pm}|=|g_n|[1+Y_n(z_n^{\pm})]^2$  (see [KK1])
implies $2|\m_n^{\pm}|\ge |g_n|$, which yields the first inequality
in \er{2.10}. Moreover, using \er{3.36} we obtain the estimate
$2|\m_n^{\pm}|\le c^{18}|g_n|$, which  gives the second one in
\er{2.10}.\ \ $\BBox$

 Recall that for a compact subset $\O\ss\C $ the analytic capacity
is given by
\[
\lb{2.11} \cC=\cC (\O)=\sup\lt[|f'(\iy)|: f \ is \ analytic\ in  \
\C \sm \O; \ \ \ |f(k)|\le 1,\ k\in\C \sm \O \rt],
\]
where $f'(\iy)=\lim_{|k|\to\iy}k(f(k)-f(\iy))$. We will use the well
known Theorem (see \cite{Iv}, \cite{Po})

\no {\bf Theorem }    {\it ( Ivanov-Pommerenke ).
 Let $E\ss\R $ be compact. Then  the analytic capacity
$\cC (E)=|E|/4$, where $|E|$ is the Lebesgue measure (the length) of
the set  $E$. Moreover, the Ahlfors function $f_E$ (the unique
function, which gives  $\sup$ in the definition of  the analytic
capacity) has the following form:}
\[
\lb{2.12}
f_E(z)=\frac{\exp{(\frac12\f_E(z))}-1}{\exp{(\frac12\f_E(z))}+1},\ \
\ \ \ \ \  \f_E(z)=\int_E\frac{dt}{z-t}; \ \ z\in\C \sm E.
\]
We will use the following simple remark: Let $S_1, S_2,\dots, S_N$
be disjoint continua in the  plane $\C ;\ \ D=\C
\sm\cup_{n=1}^NS_n$. Introduce the class $\S '(D)$ of the conformal
mapping $w$ from the domain $D$ onto $\C $ with the following
asymptotics: $w(k)=k+[Q(w)+o(1)]/k,\  k\to\iy$. If $\O\ss\C $ is
compact; $ D=\C \sm \O,\ g\in\S '(D)$, then  $\cC (\O)=\cC (\C \sm
g(D))$. It follows immediately from the definition of the analytic
capacity.

Let $\ell^2_{fin}\ss \ell^2$ be the subset of finite sequences of
non negative numbers. Then, using  the Ivanov-Pommerenke Theorem
 and the last remark we obtain
$$
\|g(h)\|_1=\cC (\G (h)),\ \ where \ \ \  h\in\ell_{fin}^2, \ \G
(h)=\cup [u_n-ih_n,u_n+ih_n];\ \ \ \
$$

Now we estimate the Dirichlet integral $I_D(h)=2Q_0(h)$ for the case
$u_*\ge 0$, using a geometric construction.

\begin{theorem}\lb{T2.4}
Let $h\in \ell^\iy, h_n\to 0 $ as $|n|\to\iy;$ and let $\wt h=\wt
h(h)$ be given by

if $h=0$, then $\wt h=0$,

if $ h\ne 0$, let an integer $n_1 $ be such that $\wt
h_{n_1}=h_{n_1}=\max_{n\in\Z }h_n>0$; assume that  the numbers
 $h_{n_1},h_{n_2},\dots ,h_{n_k} $ are defined,  then $n_{k+1} $ is given by
\[
\lb{2.15} \wt h_{n_{k+1}}=h_{n_{k+1}}=\max_{n\in B}h_n>0, \ \
B=\{n\in\Z : |u_n-u_{n_s}|>h_{n_s}, 1\le s\le k\},
\]
and let $\wt h_n=0$, if $n \notin \{n_k, k\in \Z\}$.

 Then the following estimates hold:
\[
\lb{2.16} \frac1{\pi^2}\|\wt h\|^2_2\le Q_0(h)=
{I_D(h)\/2}\le\frac{2\sqrt2}\pi\|\wt h\|^2_2.
\]
\end{theorem}

\no {\bf Proof.}  The Lindel\"of principal yields $Q_0(\wt h)\le
Q_0(h) $. On the other hand
 open squares $P_k=(u_{n_k}-t_k,u_{n_k}+t_k)\times(-t_k,t_k),
t_k\ev h_{n_k}=\wt h_{n_k}, k\in Z $, does not overlap. Then
applying \er{3.3} to the function
 $(z(k,\wt h)-k) $ and $P_k$, we obtain
\[
\lb{3.29} 2t_k^2\le\pi\int \int_{P_k}|z'(k,\wt h)-1|^2\,dudv,
\]
and
$$
Q_0(\wt h)=\frac12I_D(\wt
h)\ge\frac1{\pi^2}\sum_{k\ge1}t_k^2=\frac1{\pi^2} \|\wt h\|^2_2
$$
which yields the first estimate in \er{2.16}.

   Let $\O_k=\{n\in\Z : u_n\in[u_{n_k}-t_k,u_{n_k}+t_k]\} $.
By the Lindel\"of principal, the  gap length $|g_{n_k}|$ such that
 $[u_n,u_n+ih_n], n\in \O_k $ increases if we take off all another slits.
By the Ivanov-Pomerenke Theorem (see above), the sum of new gap
lengths equals to $4\times capacity $ of the set
 $E=\cup_{n\in \O_k}[u_n-h_n,u_n+h_n] $, which is less than
the diameter of the set $E$. Then $\sum_{n\in
\O_k}|g_n(h)|\le2\sqrt2t_k$, and using the last estimate we obtain
\[
\lb{3.31} \pi Q_0(h)\le\sum_{n\in\Z }h_n|g_n|\le
\sum_{k\ge1}\sum_{n\in \O_k}h_n|g_n| \le \sum_{k\ge1}t_k\sum_{n\in
\O_k}|g_n|\le 2\sqrt2\sum_{k\ge1}t_k^2=2\sqrt2\|\wt h\|^2_2,
\]
since $h_n\le t_k, n\in \O_k $ and the diameter of the set $E $ is
less than or equals $2\sqrt2t_k $.  $\BBox$

Note that the proved Theorem shows that estimates \er{3.7},\er{3.8}
hold true for the weaker conditions on the sequence $u_n,n\in \Z$.

\no {\bf Acknowledgments.}
E. Korotyaev was partly supported by DFG project BR691/23-1.
The various parts of this paper were written at the Mittag-Leffler Institute,
Stockholm  and in the Erwin Schr\"odinger Institute for Mathematical Physics,
Vienna, E. Korotyaev is grateful to the Institutes for the hospitality.

\end{document}